\documentclass[12pt]{article}
\usepackage{amsmath,amssymb}
\usepackage{titletoc}
\usepackage{latexsym}
\usepackage{amsmath}
\usepackage{amssymb}
\usepackage{multicol}
\usepackage{graphics}
\usepackage{graphicx}
\usepackage{tikz}
\usetikzlibrary{calc}
\usepackage{subfigure}
\usepackage{indentfirst}
\usepackage{epsfig}
\usepackage{mathrsfs}
\newcommand{\p}{\partial}

\title{An inverse problem of determining fractional orders in a fractal solute transport model
\author{Gongsheng Li\footnote{Corresponding author, Email: ligs@sdut.edu.cn} \hspace{2mm} Xianzheng Jia  \hspace{2mm} Wenyi Liu \hspace{2mm} Zhiyuan Li\\
{\small School of Mathematics and Statistics, Shandong University of Technology}\\
{\small Zibo, Shandong 255049, China}
}
\date{}
}
\begin{document}

\maketitle

\begin{center}
\begin{minipage}{13cm}
{\bf Abstract}: A fractal mobile-immobile (MIM in short) solute transport model in porous media is set forth, and an inverse problem of determining the fractional orders by the additional measurements at one interior point is investigated by Laplace transform. The unique existence of the solution to the forward problem is obtained based on the inverse Laplace transform, and the uniqueness of the inverse problem is proved in the real-space of Laplace transform by the maximum principle, and numerical inversions with noisy data are presented to demonstrate a numerical stability of the inverse problem.
\end{minipage}
\end{center}
\section{Introduction}
Solute transport in porous media is a complicated process involving in physical/chemical and biological reactions with fluid mechanics, and the traditional models are the advection-dispersion equations and the mobile-immobile (MIM in short) solute transport models. The MIM model describes the hydrodynamic behavior in the mobile zone and the mass transfer process between the mobile zone and the immobile zone, which can characterize the physical/chemical non-equilibrium of solute transport in heterogeneous porous media. Although the physical and chemical non-equilibrium models are based on different concepts, they can be described by the same mathematical equation in dimensionless form, see \cite{pang99, tor95} for instance. A MIM solute transport undergoing linear sorption without degradations and source/sink reactions is expressed via:
$$
\left\{
\begin{array}{llll}
\beta R \frac{\p C_1}{\p t}=\frac{1}{P} \frac{\p^2 C_1}{\p x^2}-\frac{\p C_1}{\p x}-\omega (C_1- C_2),\\
(1-\beta) R \frac{\p C_2}{\p t}=\omega (C_1- C_2),
\end{array}
\right.
\eqno{(1.1)}
$$
where $C_1, C_2$ are the dimensionless solute concentrations in the equilibrium and non-equilibrium sites respectively, $P>0$ is the Pelect number, and $R\geq 1$ is the retardation factor due to the sorption, and $\beta\in (0, 1)$ is a partitioning coefficient between the equilibrium and non-equilibrium phases, and $\omega>0$ is the first-order mass transfer rate. \\
\indent The system (1.1) is a classical integer-order MIM model for solute transport in porous media which has been studied and applied widely by hydrogeologists not only in laboratory but also in field tests, see \cite{ben09, gao10, lix20, luc18, van89, zhang04} for instance. However, there were some researches in the last decades constantly indicated that fractional differential equations could be more suitable than those of classical models to describe non-Darcian flow or anomalous diffusion in some special environment, especially in low-permeability porous media, see \cite{cap04, han01, mil63, obe17, rag12, zhou18} for instance. The solute mass transfer or the chemical reaction in a heavy heterogeneous porous media is not an instantaneous process but a longtime dynamical behavior due to the memory effect, in which case fractional diffusion equations incorporating with the memory effect are expected to describe the anomalous diffusion processes, see \cite{bae07, ben00, ger06, kel19, sch03, zhangy09, zhou19}, for instance.\\
\indent This paper is devoted to a modified model of (1.1) by replacing the first-order derivatives on time in the model with Caputo fractional derivatives, which results in a novel fractal MIM solute transport system. Such a fractal MIM model can describe some anomalous diffusion behaviors in the mobile zone and dynamical processes with memory effect in the immobile zone especially in low-permeability porous media. It is important to study the solution of the coupled fractal model, however, it is of the same importance to identify and determine those unknown parameters in the model with suitable additional information, which leads to corresponding inverse problems in the fractal MIM solute transport.\\
\indent On the research of the forward problem like the system (1.1) including the fractional-order forms, the method of Laplace transform is often utilized to deduce an expression of the solution in frequency domain, and then numerical solution is obtained by approximating the inverse Laplace transform (see \cite{sch03} for instance). However, there are no theoretical analysis to the solution's properties in mathematics. For the fractal MIM solute transport model in this paper, we will give the unique existence of the solution to the forward problem also by the method of Laplace transform, where a bounded estimate for the mapping function of Laplace transform is established to ensure the convergence of the contour integral.\\
\indent As for inverse problems associated with a system of fractional differential equations, there are still few studies in the known literatures. For the researches on inverse problems in one fractional diffusion equation, we refer to \cite{cheng09, lig13, liujj10, saka11, yama12} for some early work, and recently see \cite{kian18, kian20, liz20, sunc20, sunl19, xian20, yama21, zheng19} and the references therein.\\
\indent The fractional order in a fractal model is a key parameter to characterize the heavy-tail sub-diffusion of the solute with memory effect. However, it is always unknown in advance which leading to inverse problems of identifying the fractional order. We will consider an inverse problem of determining the two fractional orders in the fractal MIM solute transport system using the additional data measured at one interior point. The uniqueness of the inverse problem is proved by the Laplace transform method under the condition that infinite measurements can be obtained at the space point. Such condition seems to be unreasonable for real-life problems, but it gives us an approach to the Laplace transform for the solution on $t\geq 0$, and it only needs a set of limited data on concrete numerical inversions. Furthermore, based on the finite difference solution of the forward problem, numerical inversions with noisy data are presented by using a modified Levenberg-Marquart algorithm.\\
\indent The rest of the paper is organized as follows.\\
\indent In Section 2, some preliminaries on the Laplace transform and the maximum principle are given, and in Section 3 a fractal MIM solute transport model is introduced, and the unique solvability of the forward problem is derived based on the inverse Laplace transform. In Section 4, an inverse problem of determining the factional orders is considered, and its uniqueness is proved by the maximum principle in the real space of the Laplace transform. In Section 5, numerical inversions with noisy data are presented to demonstrate a numerical stability of the inverse problem, and concluding remarks are given in section 6.
\section{Preliminaries}
In this section we give some preliminaries on the Laplace transform and its inverse transform of a real-valued function, and the maximum principle of elliptic operator.
\subsection{Basic facts on the Laplace transform}
In this subsection, the function $f(t)$ is assumed to be the first-order differentiable on $t\in [0, \infty)$ such that the first-order derivative $f'(t)$ and the $\alpha$-order fractional derivative $\p_t^\alpha f(t)$ ($0<\alpha<1$) exist. The function $\bar{f}(s)$ of the complex variable $s$ defined by
$$
\bar{f}(s)={\cal L}\{f(t); s\}=\int_0^\infty \exp(-s t) f(t) dt,\eqno{(2.1)}
$$
is called the Laplace transform of $f(t)$ ($t\geq 0$), where $f(t)$ satisfies the growth condition $|f(t)|\leq M \exp(c_0 t)$ as $t\rightarrow \infty$, and $M, c_0$ are positive constants. \\
\indent If confining the parameter $s$ in the real space of $s>c_0$, we can get the sign of the Laplace transform function.\\
\textbf{Lemma 2.1} \quad Assume that the function $f(t)$ is nonnegative for $t\in [0, \infty)$ and satisfies the growth condition, then there holds
$$
\bar{f}(s)\geq 0, \qquad s> c_0,\eqno{(2.2)}
$$
and $\bar{f}(s)\rightarrow 0$ as $s\rightarrow +\infty$.\\
\textbf{Proof}
\quad Obviously, if $f(t)\geq 0$ and  $s> c_0>0$, there must have $\bar{f}(s)\geq 0$ by (2.1). Furthermore, there holds
$$
\bar{f}(s)\leq M \int_0^\infty \exp(-s t) \exp(c_0 t) dt=\frac{M}{s-c_0}\rightarrow 0,\ s\rightarrow+\infty.
\eqno{(2.3)}
$$
\indent The inverse Laplace transform of the function $f(t)$ is defined via:
$$
f(t)={\cal L}^{-1}\{\bar{f}(s); t\}=\int_{s_0-i\infty}^{s_0+i\infty} \bar{f}(s) \exp(s t) ds, \eqno{(2.4)}
$$
where $s_0=Re(s)>c_0$.\\
\textbf{Lemma 2.2}\quad  If the Laplace transform function $\bar{f}(s)$ satisfies the condition
$$
|\bar{f}(s)|\leq \frac{C}{|s|},\eqno{(2.5)}
$$
where $C>0$ is a constant independent of $s$ and $Re(s)>c_0$, then the contour integral in (2.4) is convergent at each given $t>0$, and the inverse Laplace transform is well-defined.\\
\textbf{Proof}\quad See the Appendix.\\
\indent Finally we give the Laplace transform of the Caputo fractional derivative $\p_t^\alpha f(t)$ ($0<\alpha<1$). The Caputo fractional derivative $\p_t^\alpha f(t)$ for $0<\alpha<1$ is defined by
$$
\p_t^\alpha f(t)=\frac{1}{\Gamma(1-\alpha)} \int_0^t \frac{f'(\tau)}{(t-\tau)^\alpha} d\tau,\eqno{(2.6)}
$$
where $\Gamma(\cdot)$ denotes the Gamma function, see \cite{kil06, pod99} for detailed definitions and properties of fractional derivatives.\\
\indent On performing Laplace transform for a fractional derivative function, some regularity is needed for the performed function, see \cite{kub20} for detailed analysis. We set
$$
W^{1,1}(0, \infty):=\{f(t)\in L^1(0,\infty); f'(t)\in L^1(0,\infty)\};\eqno{(2.7)}
$$
and for $\alpha\in (0, 1)$, we set
$$
W_\alpha(0, \infty):=\{f(t)\in W^{1,1}(0, \infty);\ t^{1-\alpha} f'(t)\in L^\infty(0, \infty)\}.\eqno{(2.8)}
$$
Next for $\alpha\in (0, 1)$, we set
$$
V_\alpha(0, \infty):=\{f(t)\in W_\alpha(0, \infty);\ \exists M, c_0>0\  \hbox{such that}\  |f(t)|\leq M e^{c_0 t}\}.\eqno{(2.9)}
$$
Now for $f(t)\in V_\alpha(0, \infty)$, we can define the Laplace transform of the fractional derivative $\p_t^\alpha f$ as follows:
$$
{\cal L}\{\p_t^\alpha f(t); s\}=\int_0^\infty e^{-s t} \p_t^\alpha f(t) dt,\ Re(s)>c_0,\eqno{(2.10)}
$$
and there holds
$$
{\cal L}\{\p_t^\alpha f(t); s\}=s^\alpha \bar{f}(s)-s^{\alpha-1} f(0), \eqno{(2.11)}
$$
where $\bar{f}(s)$ denotes the Laplace transform of $f(t)$ on $t\in [0, \infty)$.

\subsection{Maximum principle of elliptic operator}
\textbf{Lemma 2.3}(\cite{spe81})\quad Let $I$ be a bounded interval in ${\bf R}$, and $u=u(x)$ be a nonconstant solution of
$$
a(x) u''+b(x) u'+h(x) u\geq 0,\  x\in I,\eqno{(2.12)}
$$
where the coefficients $a(x), b(x)$ and $h(x)$ are bounded and $h(x)\leq 0$ in $I$, and there exists a constant $a_0>0$ such that $a(x)\geq a_0>0$ in $I$. Then a nonnegative maximum of $u$  can only occur on $\p I$, and $du/d\nu>0$ there, where $\nu$ denotes a normal vector pointing outward at the boundary.\\
\textbf{Corollary 2.1}\quad Let $I=(0, 1)$. Under the conditions of Lemma 2.3, suppose further that $u(0)=0$ and $u'(1)=0$, then there must have $u(x)\leq 0$, $x\in I$.\\
\textbf{Proof}\quad By Lemma 2.3, the solution $u$ can not attain its maximum at $x=1$ since $u'(1)=0$, and it has to assume its maximum at
$x=0$, i.e., there is $u(x)\leq 0$ in $I$.
\section{The fractal MIM model}
\subsection{The forward problem}
Consider the solute transport model (1.1) in a 1D finite space domain but in the infinite time domain. Let $\Omega=(0, 1)$  by dimensionless and $\Omega_\infty=\Omega\times (0, \infty)$. Assume that the solute transport and diffusion begins in the mobile phase, and the solute variation in the immobile is a dynamical process due to the low-permeability and heavy heterogeneity of the porous media. Then it could be more suitable that the solute diffusion in the immobile zone is described by a time-fractional differential equation. Correspondingly, the advection-diffusion processes in the mobile zone can also be governed by a time-fractional advection-diffusion equation. In addition, assume that there are the first-order degrading reactions in the two zones respectively, and a fractal MIM model for reactive solute transport for $(x, t)\in \Omega_\infty$ is established as follows on the basis of (1.1):
$$
\left\{
\begin{array}{lll}
\beta R_1 \p_t^\alpha u_1=\frac{1}{P} \frac{\p^{2} u_1}{\p x^{2}}-\frac{\p u_1}{\p x}-\omega (u_1-u_2)-\lambda u_1,\\
(1-\beta) R_2 \p_t^\gamma u_2 = \omega (u_1-u_2)-\mu u_2,
\end{array}
\right.
\eqno{(3.1)}
$$
where $u_1=u_1(x,t)$ and $u_2=u_2(x,t)$ denote the solute concentrations in the mobile and the immobile zones respectively, and $R_1, R_2>1$ are the retardation coefficients with time-scale actions, and $\lambda, \mu>0$ are the first-order degradation coefficients (or the coefficients of zero-order derivatives in mathematics) in the mobile and immobile zones respectively; $\beta\in (0, 1)$ and $\omega>0$ are the same meanings as denoted in (1.1), and $\p_t^\alpha u_1$ ($0<\alpha<1$) and $\p_t^\gamma u_2$ ($0<\gamma<1$) denote the Caputo fractional derivatives on time $t>0$. Here the fractional orders $\alpha$ and $\gamma$ are the indexes describing the sub-diffusion characteristics with long-time memory in the mobile and immobile zones, respectively. \\
\indent For the model (3.1), the initial condition is given as:
$$
u_1(x, 0)=0,\quad  u_2(x, 0)=0,\quad 0\leq x\leq 1, \eqno{(3.2)}
$$
which means that the concentration of the solute in the studied region is zero at the initial stage. The boundary condition at $x=0$ is given as
$$
u_1(0, t)=1,\quad u_2(0, t)=0, \quad 0< t< \infty, \eqno{(3.3)}
$$
which implies that the left-hand side of the region in the mobile is an input source for $t>0$. The boundary condition at $x=1$ is impermeable, which is given by
$$
\frac{\p u_1}{\p x} (1, t)=0,\quad  \frac{\p u_2}{\p x} (1, t)=0, \quad 0< t< \infty. \eqno{(3.4)}
$$
\indent As a result, we get a coupled system composed by the fractal MIM solute transport model (3.1) with the initial boundary value conditions (3.2)-(3.4), which is called the forward problem. We consider the unique solvability of the forward problem by Laplace transform method in the next subsection.
\subsection{Existence of the solution}
Due to physical reasons, the only solutions of (3.1) we are interested in are the bounded and nonnegative ones in $\Omega_\infty$, and according to the background of solute transport in  porous media and the physical/chemical laws, the parameters in the model (3.1) satisfy the following natural condition throughout this paper:\\
$(A1)\quad 0<\alpha,\ \gamma<1,\ 0<\beta<1,\  R_1, R_2\geq 1,\ P>0,\ \omega>0,\ \lambda>0,  \mu>0.$\\
\indent Suppose that the Laplace transforms on $t\geq 0$ for all functions in the system (3.1) are existed. By performing Laplace transform for the system (3.1), and thanks to the formula (2.11) and the homogeneous initial condition (3.2),  we get
$$
\left\{
\begin{array}{lll}
\beta R_1 s^\alpha \bar{u}_1=\frac{1}{P} \frac{d^{2} \bar{u}_1}{d x^{2}}-\frac{d \bar{u}_1}{d x}-\omega (\bar{u}_1-\bar{u}_2)-\lambda \bar{u}_1,\\
(1-\beta) R_2 s^\gamma \bar{u}_2 = \omega (\bar{u}_1-\bar{u}_2)-\mu \bar{u}_2,
\end{array}
\right.
\eqno{(3.5)}
$$
where $\hbox{Re}(s)>0$ due to the boundedness of the solution. Since there is
$$
\bar{u}_2=\frac{\omega}{(1-\beta) R_2 s^\gamma+ \omega+\mu} \bar{u}_1,\eqno{(3.6)}
$$
we have
$$
a \frac{d^{2} \bar{u}_1}{d x^{2}}-\frac{d \bar{u}_1}{d x}+ b \bar{u}_1= 0,\eqno{(3.7)}
$$
where $a=\frac{1}{P}$, and
$$
b=-\beta R_1 s^\alpha-\omega-\lambda+\frac{\omega^2}{(1-\beta)R_2 s^\gamma+\omega+\mu}.\eqno{(3.8)}
$$
\indent It is noted that Eq.(3.7) is the second-order ordinary differential equation on $x\in \Omega$ with constants coefficients. By (3.3) and (3.4) the boundary conditions are given as
$$
\bar{u}_1(0;s)=\frac{1}{s},\   \bar{u}'_1(1;s)=0.\eqno{(3.9)}
$$
In the follows we give a solution's representation for the boundary value problem (3.7), (3.9) by the eigenvalue method.\\
\indent By using the trigonometric representation of complex number, and noting $\hbox{Re}(s)>0$, there must have $\hbox{Re}(b)<0$ by (3.8).
By solving the characteristic equation ($a>0$)
$$
a \eta^2-\eta+ b=0,
$$
we get
$$
\eta_{1}=\frac{1+ \sqrt{1-4ab}}{2a},\  \eta_{2}=\frac{1- \sqrt{1-4ab}}{2a},\eqno{(3.10)}
$$
where $\hbox{Re}(\eta_1)>0$ and $\hbox{Re}(\eta_2)<0$ due to $\hbox{Re}(b)<0$ and $a>0$. So the solution of the problem (3.7), (3.9) is expressed by
$$
\bar{u}_1(x;s)=c_1 e^{\eta_1 x}+ c_2 e^{\eta_2 x},\eqno{(3.11)}
$$
where
$$
c_1=\frac{\eta_2 s^{-1}}{\eta_2-\eta_1 e^{\eta_1-\eta_2}}; c_2=\frac{\eta_1 s^{-1}}{\eta_1-\eta_2 e^{\eta_2-\eta_1}}. \eqno{(3.12)}
$$
Together with (3.6) follows the expression of $\bar{u}_2(x;s)$. In order to utilize the inverse Laplace transform to obtain the solution of the forward problem, we need a bounded estimate for $\bar{u}_1$ given by (3.11).\\
\textbf{Lemma 3.1}\quad For the mapping function of Laplace transform given by (3.11), there holds
$$
|\bar{u}_1| \leq \frac{1}{\left|s\right|} C,  \eqno{(3.13)}
$$
where $C>0$ is a constant independent of $s$, and $\hbox{Re}(s)>0$.\\
\textbf{Proof} \quad As indicated in the above, there is $\hbox{Re}(b)<0$ for $\hbox{Re}(s)>0$. From (3.11) we have
$$
\begin{aligned}
|\bar{u}_1| & \leq \frac{1}{|s|} \left|\frac{\eta_{2} e^{\eta_{2}}}{\eta_{2} e^{\eta_{2}}-\eta_{1} e^{\eta_{1}}}\right|\left|e^{\eta_{1} x}\right|+\frac{1}{|s|}\left|\frac{\eta_{1} e^{\eta_{1}}}{\eta_{1} e^{\eta_{1}}-\eta_{2} e^{\eta_{2}}}\right|\left| e^{\eta_{2} x}\right| \\
& =\frac{1}{|s|}(I_1+I_2),
\end{aligned}
\eqno{(3.14)}
$$
where
$$
I_1=\frac{1}{\left|1-\frac{\eta_{1}}{\eta_{2}}e^{\eta_{1}-\eta_{2}}\right|}\left|e^{\eta_{1}x}\right|;\
I_2=\frac{1}{\left|1-\frac{\eta_{2}}{\eta_{1}} e^{\eta_{2}-\eta_{1}}\right|}\left|e^{\eta_{2} x}\right|.
$$
\indent For the estimates of $I_1$ and $I_2$, we need the properties of $\eta_{1}$ and $\eta_{2}$. By (3.10) there are
$$
\frac{\eta_{1}}{\eta_{2}}=\frac{(1+\sqrt{1-4ab})^{2}}{4ab},\  \frac{\eta_{2}}{\eta_{1}}=\frac{(1-\sqrt{1-4ab})^{2}}{4ab};\eqno{(3.15)}
$$
and
$$
\eta_{1}-\eta_{2}=\frac{\sqrt{1-4ab}}{a},\ \eta_{1}+\eta_{2}=\frac{1}{a}.\eqno{(3.16)}
$$
By the expression of $b$ given in (3.8),there holds
$$
b\rightarrow -\omega-\lambda+\frac{\omega^2}{\omega+\mu},\ \hbox{as}\ s\rightarrow 0,
$$
which means that the norm of the coefficient $b$ is lower bounded.  Now we estimate the term of $I_1$.\\
\indent Firstly by (3.15) and noting $\hbox{Re}(b)<0$, there exists a positive constant $C_1$ such that $|\frac{\eta_{1}}{\eta_{2}}|=1+C_1>1$, and there is $\hbox{Re}(\eta_1-\eta_2)>0$.
Then there holds
$$
|\frac{\eta_{1}}{\eta_{2}} e^{\eta_{1}-\eta_{2}}|=|\frac{\eta_{1}}{\eta_{2}}|\cdot |e^{\eta_{1}-\eta_{2}}|=(1+C_1) e^{\hbox{Re}(\eta_1-\eta_2)}>1.\eqno{(3.17)}
$$
Next by (3.16) there is
$$
\eta_1=\frac{1}{2a}+\frac{\eta_1-\eta_2}{2}.
$$
Noting $x\in [0, 1]$, we have
$$
|e^{\eta_1} x|=|e^{\frac{x}{2a}}| |e^{\frac{(\eta_1-\eta_2)x}{2}}|\leq e^{\frac{1}{2a}} e^{\frac{\hbox{Re}(\eta_1-\eta_2)}{2}}.
$$
Therefor we get
$$
I_1\leq \frac{|e^{\eta_1 x}|}{\left|\frac{\eta_{1}}{\eta_{2}} e^{\eta_{1}-\eta_{2}}\right|-1}
\leq e^{\frac{1}{2a}} \frac{e^{\frac{\hbox{Re}(\eta_1-\eta_2)}{2}}}{(1+C_1) e^{\hbox{Re}(\eta_1-\eta_2)}-1}
\leq e^{\frac{1}{2a}} \frac{e^{\hbox{Re}(\eta_1-\eta_2)}}{e^{\hbox{Re}(\eta_1-\eta_2)}-1},
\eqno{(3.18)}
$$
which implies that there exists a constant $C_2>0$ such that $I_1\leq C_2$.\\
\indent Similarly there exists a constant $C_3>0$ such that $I_2\leq C_3$, thus the assertion (3.13) is valid, and the proof is completed.\\
\indent With the above lemma, we are ready to give the unique existence of the solution to the forward problem.\\
\noindent\textbf{Theorem 3.1}\quad The forward problem (3.1) with (3.2)-(3.4) has a unique solution in $\Omega_{\infty}$.\\
\textbf{Proof}\quad We only need to prove the existence of the inverse Laplace transform on $\bar{u}_1$. By Lemma 3.1, there is $|\bar{u}_1(x, s)|  \leq \frac{1}{\left|s\right|} C$ for $\hbox{Re}(s)>0$. Therefore utilizing Lemma 2.2, the contour integral
$$
\frac{1}{2\pi i} \int_{s_{0}-i\infty}^{s_{0}+i\infty} \bar{u}_1(x, s) e^{s t} ds,
$$
is convergent for $(x, t)\in \Omega_{\infty}$, which is the solution $u_1(x, t)$, i.e., there is
$$
u_1(x, t)=\frac{1}{2\pi i} \int_{s_{0}-i\infty}^{s_{0}+i\infty} \bar{u}_1(x, s) e^{s t} ds, \eqno{(3.19)}
$$
where $s_0=\hbox{Re}(s)>0$. Similarly we can get the expression of the solution $u_2(x, t)$.\\
\indent This theorem gives the unique existence of the solution to the forward problem, however, the space for the solution is not deduced since the solution's regularity is still open. It is noted that the solution in the form of Laplace transform is not convenient to practice and application due to the expensive computational cost of the contour integral. Nevertheless, it is meaningful in mathematics we get the existence of the solution to the forward problem, and we we will give a finite difference solution in Section 5.\\
\indent In the follows, we consider an inverse problem of identifying the fractional orders $\alpha\in (0, 1)$ and $\gamma\in (0, 1)$ by the additional measurements on $u_1$ at one interior point, and we will prove its uniqueness also based on the Laplace transform, and perform numerical inversions by the Levenberg-Marquart algorithm together with homotopy technique.
\section{The inverse problem of fractional orders}
\subsection{The inverse problem}
When the model (3.1) is applied to study a real solute transport problem in a heterogeneous porous media, some model parameters are always unknown, such as the fractional order, the mass transfer rate, etc. Suppose that the fractional orders $\alpha$ and $\gamma$ are unknown, we are to determine them by some additional measurements at one interior point in the mobile zone.\\
\indent The additional condition is given as
$$
u_1(x_0, t),\ t>0,\eqno{(4.1)}
$$
where $x_0\in \Omega$ is a fixed point. Based on the above discussions, the inverse problem is to identify the two fractional orders $\alpha\in (0, 1)$ and $\gamma\in (0, 1)$ using the overposed condition (4.1) based on the forward problem (3.1), (3.2)-(3.4).
\subsection{The uniqueness}
The uniqueness of an inverse problem is important not only for theoretical analysis but also for numerical computations. We will prove the uniqueness in the mapping space of Laplace transform using the maximum principle of elliptic operator. \\
\indent An inverse problem is often investigated in an admissible set of the unknowns. For the considered inverse problem, we assume that the unknown parameters satisfy the natural condition $(\alpha, \gamma)\in S_{ad}$, where $S_{ad}$ is given by
$$
S_{a d}=\{(\alpha, \gamma): 0<\alpha<1, 0<\gamma<1\}.\eqno{(4.2)}
$$
For any given $(\alpha, \gamma)\in S_{a d}$, denote ${u}_{1}^{\alpha, \gamma}(x, t)$ as the solution of the forward problem in the mobile zone, and $u_{2}^{\alpha,\gamma}(x,t)$ the solution in the immobile zone. The solution should have some regularity so as to ensure to perform Laplace transforms for the solution itself and its derivatives, including the fractional-order derivatives. In the real-space of the Laplace transform, we can prove a uniqueness for the inverse fractional order problem.\\
\textbf{Theorem 4.1}\quad Assume that $u_{1}^{\alpha,\gamma}, u_{2}^{\alpha,\gamma}$ are the nonnegative and bounuded solutions of the forward problem corresponding to $(\alpha, \gamma)\in S_{ad}$ in the mobile and immobile zones respectively, and $x_{0}\in \Omega$ is a fixed interior point.
If $u_{1}^{\alpha_1,\gamma_1}(x_0, t)=u_{1}^{\alpha_2,\gamma_2}(x_{0}, t)$ for $t>0$ and $(\alpha_i, \gamma_i)\in S_{ad}$ ($i=1, 2$), then there holds
$\alpha_1=\alpha_2$ and $\gamma_1=\gamma_2$.\\
\textbf{Proof}
\quad  By utilizing the Laplace transform and noting the homogeneous initial condition, there hold
$$
\beta R_1 (s^\alpha \hat{u_1}^{\alpha,\gamma})=\frac{1}{P} \frac{d^{2} \hat{u_1}^{\alpha,\gamma}}{d x^{2}}-\frac{d \hat{u_1}^{\alpha,\gamma}}{d x}-\omega(\hat{u_1}^{\alpha,\gamma}-\hat{u_2}^{\alpha,\gamma})-\lambda \hat{u_1}^{\alpha,\gamma},
\eqno{(4.3)}
$$
and
$$
(1-\beta) R_2 (s^{\gamma} \hat{u_2}^{\alpha,\gamma})=\omega(\hat{u_1}^{\alpha,\gamma}-\hat{u_2}^{\alpha,\gamma})-\mu \hat{u_2}^{\alpha,\gamma}.
\eqno{(4.4)}
$$
From (4.4) there is
$$
\hat{u_2}^{\alpha,\gamma}=\frac{\omega \hat{u_1}^{\alpha,\gamma}}{(1-\beta) R_2 \  s^{\gamma}+\omega+\mu}.\eqno{(4.5)}
$$
Substituting (4.5) into (4.3) we get
$$
\frac{1}{P} \frac{d^{2} \hat{u_1}^{\alpha,\gamma}}{d x^{2}}-\frac{d \hat{u_1}^{\alpha,\gamma}}{d x}+\left\{\frac{\omega^{2}}{(1-\beta) R_2 s^{\gamma}+\omega+\mu}-\omega-\lambda-\beta R_1 s^\alpha\right\} \hat{u_1}^{\alpha,\gamma}=0.
\eqno{(4.6)}
$$
\indent Now for $(\alpha_i, \gamma_i)\in S_{ad}$ ($i=1, 2$), there hold the expressions for $\hat{u_1}^{\alpha_1,\gamma_1}$ and $\hat{u_1}^{\alpha_2, \gamma_2}$ corresponding to (4.6). Let $U(x)=\hat{u_1}^{\alpha_1,\gamma_1}-\hat{u_1}^{\alpha_2,\gamma_2}$ and assume that $\alpha_1>\alpha_2$. There holds for $x\in \Omega$
$$
\frac{1}{P} \frac{d^{2} U}{d x^{2}}-\frac{d U}{d x}+ c U =d,\eqno{(4.7)}
$$
where
$$
c=-\omega-\lambda-\beta R_1 s^{\alpha_1} +\frac{\omega^{2}}{(1-\beta) R_2 s^{\gamma_1}+\omega+\mu},\eqno{(4.8)}
$$
and
$$
d=\hat{u_1}^{\alpha_2,\gamma_2}\left\{\beta R_1 (s^{\alpha_1}-s^{\alpha_2})+\frac{\omega^{2} (1-\beta) R_2 (s^{\gamma_1}-s^{\gamma_2})}{[(1-\beta) R_2 s^{\gamma_{1}}+\omega+\mu][(1-\beta) R_2 s^{\gamma_{2}}+\omega+\mu]}\right\},
\eqno{(4.9)}
$$
and the boundary conditions are $U(0)=0$ and $U'(1)=0$.\\
\indent Let us consider the case of $s\geq s_0>0$. Thanks to the priori conditions of $\omega>0, \lambda>0, \mu>0$ and $R_1, R_2\geq 1, 1-\beta>0$, there holds
$$
\begin{array}{lll}
c&=&-\omega-\lambda-\beta R_1 s^\alpha+\frac{\omega^{2}}{(1-\beta) R_2 s^{\gamma_{1}}+\omega+\mu}\\
&\leq & -\omega-\lambda-\beta R_1 s^\alpha+\omega\\
&\leq &-\lambda-\beta R_1 s_0^\alpha<0, \quad  s\geq s_0.
\end{array}
$$
Rewrite (4.9) as
$$
d=s \hat{u_1}^{\alpha_2,\gamma_2}\left\{\beta R_1 \frac{s^{\alpha_1}-s^{\alpha_2}}{s}+\frac{\omega^{2} (1-\beta) R_2 (s^{\gamma_1}-s^{\gamma_2})}{s [(1-\beta) R_2 s^{\gamma_{1}}+\omega+\mu][(1-\beta) R_2 s^{\gamma_{2}}+\omega+\mu]}\right\}.
\eqno{(4.10)}
$$
\indent By Lemma 2.1, thanks to the nonnegative property of the solution $u_1(x, t)$, follows that $s\hat{u_1}^{\alpha_2,\gamma_2}\geq 0$ for $s\geq s_0>0$.\\
\indent By the assumption $\alpha_1>\alpha_2$ there holds $s^{\alpha_1}-s^{\alpha_2}>0$ ($s>1$) and
$$
\beta R_1 \frac{s^{\alpha_1}-s^{\alpha_2}}{s}\sim s^{-1+\alpha_1},\  s\rightarrow \infty,\eqno{(4.11)}
$$
here and in the follows, the symbol $\sim$ denotes an equivalence, $A\sim B$ means that $A/B\rightarrow \hbox{constant}$.
By the a priori conditions for the known parameters we have
$$
\frac{\omega^{2} (1-\beta) R_2(s^{\gamma_1}-s^{\gamma_2})}{s [(1-\beta) R_2 s^{\gamma_{1}}+\omega+\mu][(1-\beta) R_2 s^{\gamma_{2}}+\omega+\mu]}
\sim \frac{\omega^2}{(1-\beta) R_2} \frac{s^{-\gamma_2}-s^{-\gamma_1}}{s}, s\rightarrow \infty.\eqno{(4.12)}
$$
Since  $\gamma_1, \gamma_2\in (0, 1)$, and $s^{-\gamma_2}-s^{-\gamma_1}\rightarrow 0$ as $s\rightarrow \infty$, there holds
$$
\frac{\omega^2}{(1-\beta) R_2} \frac{s^{-\gamma_2}-s^{-\gamma_1}}{s}\sim s^{-1-\gamma^*}, s\rightarrow \infty,\eqno{(4.13)}
$$
where $\gamma^*=\min\{\gamma_1, \gamma_2\}$. Noting that
$$
\frac{s^{-1+\alpha_2}}{s^{-1-\gamma^*}}=s^{\alpha_2+\gamma^*}\rightarrow \infty, s\rightarrow \infty,\eqno{(4.14)}
$$
we get by (4.11) and (4.13)
$$
\beta R_1 \frac{s^{\alpha_1}-s^{\alpha_2}}{s}+\frac{\omega^{2} (1-\beta) R_2 (s^{\gamma_1}-s^{\gamma_2})}{s [(1-\beta) R_2  s^{\gamma_{1}}+\omega+\mu][(1-\beta) R_2 s^{\gamma_{2}}+\omega+\mu]}\geq 0,\  s\geq s_0.
\eqno{(4.15)}
$$
Together with (4.10) concludes that $d\geq 0$ for $s\geq s_0$. As a result by applying Lemma 2.3 and Corollary 2.1 to the equation (4.7) with $U(0)=0, U'(1)=0$, there holds $U(x)<0$ for $x\in \Omega$, $s\geq s_0$, and then we get
$$
U(x_0)<0,\   s\geq s_0.\eqno{(4.16)}
$$
\indent On the other hand, by the additional condition $u_1^{\alpha_1,\gamma_1}(x_0, t)=u_1^{\alpha_2,\gamma_2}(x_0, t)$ ($t>0$), we have by Laplace transform
$$
U(x_0)=\hat{u_1}^{\alpha_1,\gamma_1}(x_0,s)-\hat{u_1}^{\alpha_2,\gamma_2}(x_0, s)= 0.
\eqno{(4.17)}
$$
This is a contradiction with (4.16) and there must have $\alpha_1\leq \alpha_2$. Similarly, $\alpha_1<\alpha_2$ is impossible. Therefore $\alpha_{1}=\alpha_{2}$.\\
\indent Furthermore, denote $\alpha_1=\alpha_2:=\alpha$, we can prove $\gamma_1=\gamma_2$ by the similar arguments. Let $V(x)=\hat{u_1}^{\alpha,\gamma_1}-\hat{u_1}^{\alpha,\gamma_2}$ and assume that $\gamma_1>\gamma_2$. There holds for $x\in \Omega$
$$
\frac{1}{P} \frac{d^{2} V}{d x^{2}}-\frac{d V}{d x}+ \bar{c} V =\bar{d},\eqno{(4.18)}
$$
where
$$
\bar{c}=-\omega-\lambda +\frac{\omega^{2}}{(1-\beta) R_2 s^{\gamma_1}+\omega+\mu},\eqno{(4.19)}
$$
and
$$
\bar{d}=\hat{u_1}^{\alpha,\gamma_2}{\frac{\omega^{2} (1-\beta) R_2 (s^{\gamma_1}-s^{\gamma_2})}{[(1-\beta) R_2 s^{\gamma_{1}}+\omega+\mu] [(1-\beta) R_2 s^{\gamma_{2}}+\omega+\mu]}},
\eqno{(4.20)}
$$
and the boundary conditions are $V(0)=0$ and $V'(1)=0$.\\
\indent Also consider the case of $s\geq s_0>0$. Obviously there is $\bar{c}\leq -\lambda<0$. By the assumption of $\gamma_1>\gamma_2$ and $\gamma_1, \gamma_2\in (0, 1)$, we have
$$
s^{\gamma_1}-s^{\gamma_2}>0,\  s>1.\eqno{(4.21)}
$$
Then there holds $\bar{d}\geq 0$ for $s\geq s_0>0$, and there must have $V(x)<0$ for $x\in \Omega$ and $s\geq s_0$ also by Corollary 2.1, which leads to a contradiction with the additional condition. Thus the assumption $\gamma_1>\gamma_2$ is not valid, and similarly $\gamma_1<\gamma_2$ is not valid too. So there must have $\gamma_1=\gamma_2$. The proof is over.
\section{Numerical inversions}
This section is devoted to numerical inversions for the inverse fractional-order problem by utilizing a modified Levenberg-Marquart algorithm. On the concrete numerical computations, we only need a series of additional measurements at a limited time interval. So we can deal with the forward problem for $(x, t)\in (0, 1)\times (0, T)$, where $T>0$ is a finite number, and the additional condition is given as $\{u(x_{0}, t)\}_{0<t\leq T}$, here $x_0\in (0, 1)$ also denotes a fixed space point. For utilization of the inversion algorithm we need numerical solution of the forward problem.  Recently in \cite{liuw21}, the authors gave an implicit finite difference scheme to the forward problem, and proved its convergence and stability. For completeness of this paper, we introduce the difference scheme in the follows.
\subsection{The finite difference scheme}
Let $m, n$ be positive integers, and $h = 1/m, \tau = T/n$ be grid steps to discretize the
domain. Denote $x_{i}=i h(i=0, \cdots, m), t_{k}=k \tau(k=0, \cdots, n)$ as the grid points, and
$u_1^{i, k} \approx u_1(x_{i}, t_{k}), u_2^{i, k} \approx u_2(x_{i}, t_{k})$ as the approximations.
By the general finite difference method as used to fractional diffusion equations (see \cite{lig16, liu07, meer04} for instance), we have
\begin{equation}
\begin{array}{lll}
&\frac{\beta R_{1}}{\tau^{\alpha} \Gamma(2-\alpha)} \sum_{j=0}^{k} [u_1^{i, j+1}-u_1^{i, j}][(k+1-j)^{1-\alpha}-(k-j)^{1-\alpha}]\\
&=\frac{1}{P} \frac{u_1^{i+1, k+1}-2 u_1^{i, k+1}+u_1^{i-1, k+1}}{h^{2}}-\frac{u_1^{i, k+1}-u_1^{i-1, k+1}}{h}\\
&-\omega (u_1^{i, k+1}-\frac{u_2^{i-1, k+1}+u_2^{i+1, k+1}}{2})-\lambda u_1^{i, k+1},
\end{array} \tag{5.1}
\end{equation}
and
\begin{equation}
\begin{array}{lll}
\frac{(1-\beta) R_{2}}{\tau^{\gamma} \Gamma(2-\gamma)} \sum_{j=0}^{k}[u_2^{i, j+1}-u_2^{i, j}][(k+1-j)^{1-\gamma}-(k-j)^{1-\gamma}]\\
=\omega (\frac{u_1^{i-1, k+1}+u_1^{i+1, k+1}}{2}-u_2^{i, k+1})-\mu u_2^{i, k+1}.
\end{array} \tag{5.2}
\end{equation}
We denote $r_1=\frac{\tau^{\alpha} \Gamma(2-\alpha)}{P \beta R_{1} h^{2}}$, $r_2=\frac{\tau^{\gamma} \Gamma(2-\gamma)}{(1-\beta) R_2}$, and
\begin{equation}
\left\{\begin{array}{lll}
A=\frac{\tau^{\alpha} \Gamma(2-\alpha)}{\beta R_{1} h}+r_1,\  D=\frac{\omega \tau^{\alpha} \Gamma(2-\alpha)}{2 \beta R_1},\
E=\frac{r_2 \omega}{2}; \\
B=1+A+r_1+2 D+\frac{\tau^\alpha\Gamma(2-\alpha)}{\beta R_1}\lambda,\  F=1+2 E+r_2 \mu.
\end{array}\right. \tag{5.3}
\end{equation}
We get an implicit difference equations given as
\begin{equation}
\left\{\begin{array}{lll}
-A u_1^{i-1, k+1}+B u_1^{i, k+1}-r_1 u_1^{i+1, k+1}-D u_2^{i-1, k+1}-D u_2^{i+1, k+1} \\
=u_1^{i, k}-\sum_{j=0}^{k-1}(u_1^{i, j+1}-u_1^{i, j})[(k+1-j)^{1-\alpha}-(k-j)^{1-\alpha}], \\
-E u_1^{i-1, k+1}-E u_1^{i+1, k+1}+F u_2^{i, k+1} \\
=u_2^{i, k}-\sum_{j=0}^{k-1}(u_2^{i, j+1}-u_2^{i, j})[(k+1-j)^{1-\gamma}-(k-j)^{1-\gamma}].
\end{array}\right. \tag{5.4}
\end{equation}
Denote a new variable by
$$
U^{k}=(u_1^{1, k}, u_1^{2, k}, \cdots, u_1^{m-1, k}, u_2^{1, k}, u_2^{2, k}, \cdots, u_2^{m-1, k})^{T}, \quad k=1, 2, \cdots, n,
$$
and the initial boundary value conditions are discretized as
\begin{equation}
\begin{aligned}
U^{(0)} &=(u_1^{1,0}, u_1^{2,0}, \cdots, u_1^{m-1,0}; u_2^{1, 0}, u_2^{2, 0}, \cdots, u_2^{m-1, 0})^{T}\\
&=(0,0, \cdots, 0; 0,0, \cdots, 0)^{T},
\end{aligned} \nonumber
\end{equation}
and
\begin{equation}
\begin{array}{c}
u_1^{0, k}=1, \quad u_2^{0, k}=0, \quad k=0,1, \cdots, n; \\
u_1^{m-1, k}=u_1^{m, k}, \quad u_2^{m-1, k}=u_2^{m, k}, \quad k=0, 1, \cdots, n.
\end{array} \nonumber
\end{equation}
By rearranging (5.4) we get the difference scheme in the matrix form:
\begin{equation}
\left\{\begin{array}{l}
M U^{1}=U^{0}, \\
M U^{k+1}=N U^{k}+\sum_{j=1}^{k-1} \Psi_{j}^{k} U^{j}+N_{0} U^{0}, k=1,2, \cdots, n-1,
\end{array}\right.
\tag{5.5}
\end{equation}
where the coefficient matrix $M$ is a $2(m-1)$-order matrix defined by
\begin{equation}
M=\left(\begin{array}{ll}
M_{11} & M_{12} \\
M_{21} & M_{22}
\end{array}\right), \tag{5.6}
\end{equation}
where $M_{11}, M_{12}, M_{21}$ and $M_{22}$ are all $m-1$-order matrices given by
$$
M_{11}=\left(\begin{array}{ccccc}
{B} & {-r_1} & {0} & {\cdots} & {0}\\
{-A} & {B} & {-r_1} & {\cdots} & {0} \\
{\vdots} & {\ddots} & {\ddots} & {\ddots} & {\vdots} \\
{0} & {\cdots} & {-A} & {B} & {-r_1} \\
{0} & {\cdots} & {0} & {-A} & {B-r_1}
\end{array}\right),
M_{12}=\left(\begin{array}{ccccc}
{0} & {-D} & {0} & {\cdots} & {0} \\
{-D} & {0} & {-D} & {\cdots} & {0} \\ {\vdots} & {\ddots} & {\ddots} & {\ddots} & {\vdots}\\
{0} & {\cdots} & {-D} & {0} & {-D} \\ {0} & {\cdots} & {0} & {-D} & {-D}
\end{array}\right),
$$
$$
M_{21}=\left(\begin{array}{ccccc}
{0} & {-E} & {0} & {\cdots} & {0} \\ {-E} & {0} & {-E} & {\cdots} & {0} \\ {\vdots} & {\ddots} & {\ddots} & {\ddots} & {\vdots} \\
{0} & {\cdots} & {-E} & {0} & {-E} \\ {0} & {\cdots} & {0} & {-E} & {-E}\end{array}\right),
M_{22}=\left(\begin{array}{ccccc}
{F} & {0} & {0} & {\cdots} & {0} \\ {0} & {F} & {0} & {\cdots} & {0} \\ {\vdots} & {\ddots} & {\ddots} & {\ddots} & {\vdots} \\
{0} & {\cdots} & {0} & {F} & {0} \\ {0} & {\cdots} & {0} & {0} & {F}\end{array}\right).
$$
And the matrices $N$ and $N_{0}$ in (5.5) are all $2(m-1)$-order defined by
\begin{equation}
N=\left(\begin{array}{cc}
\left(2-2^{1-\alpha}\right) \mathbf{I} & \mathbf{O} \\
\mathbf{O} & \left(2-2^{1-\gamma}\right) \mathbf{I}
\end{array}\right),
N_{0}=\left(\begin{array}{cc}
\xi_{k}\mathbf{I} & \mathbf{O} \\
\mathbf{O} & \zeta_{k} \mathbf{I}
\end{array}\right), \tag{5.7}
\end{equation}
where $\mathbf{I}$ is the $m-1$-order identity matrix, $\mathbf{O}$
denotes the $m-1$-order zero matrix, and
$$
\xi_{k}=(k+1)^{1-\alpha}-k^{1-\alpha}, k=1, \cdots, n-1,
$$
and
$$
\zeta_{k}=(k+1)^{1-\gamma}-k^{1-\gamma}, k=1, \cdots, n-1.
$$
And the matrix $\Psi_{j}^{k}$ is defined by
\begin{equation}
\Psi_{j}^{k}=\left(\begin{array}{ll}
b_{1, j}^{k}\mathbf{I} & \mathbf{O} \\
\mathbf{O} & b_{2, j}^{k}\mathbf{I}
\end{array}\right), \tag{5.8}
\end{equation}
where
$$
b_{1,j}^{k}=2(k+1-j)^{1-\alpha}-(k-j)^{1-\alpha}-(k-j+2)^{1-\alpha},
$$
and
$$
b_{2,j}^{k}=2(k+1-j)^{1-\gamma}-(k-j)^{1-\gamma}-(k-j+2)^{1-\gamma},
$$
for $j=1, \cdots, k-1$ and $k=2, \cdots, n-1$.\\
\indent It is noted that under the natural condition (A1) given in Subsection 3.2, the coefficient matrix $M$ given by (5.6) is strictly diagonal dominant, and the finite difference scheme (5.5) is uniquely solvable. By solving the difference scheme (5.5), numerical solution of the forward problem is solved with which the modified Levenberg-Marquart algorithm is applied to give numerical inversions for the inverse problem.
\subsection{Numerical inversions}
For convenience of writing, we set $z:=(\alpha, \gamma) \in S_{ad}$ as the exact solution to the inverse problem, and $S_{ad}$ is given by (4.2), and we write the solution of the forward problem in the mobile zone as $u_1[z]$ to emphasize its dependence upon the unknown $z=(\alpha, \gamma)$. By $u_1^{\delta}(x_{0}, t)$ we denote the noisy observation data given as
$$
u_1^{\delta}(x_{0}, t)=u_1(x_{0}, t)+ \theta \delta, t\in (0, T_1], \eqno{(5.9)}
$$
where $\delta>0$ denotes the noise level, and $\theta$ is a random vector distributed in $[-1, 1]$.\\
\indent Based on the Levenberg-Marquart method, consider the following minimization problem combining with the homotopy idea:
$$
\min\limits_{z\in S_{ad}} \{(1-\kappa)\|u_1[z](x_0,t)-u_1^{\delta}(x_0, t)\|_2^2+\kappa\|z\|_2^2\},\eqno{(5.10)}
$$
where $\kappa\in (0,  1)$ is the homotopy parameter which decreases continuously from $1$ to $0$. By discretization for (5.10), and by linearization as done in the Levenberg-Marquart method, we can get a normal equation on the perturbation $\delta z$ for given $z\in S_{ad}$
$$
((1-\kappa) G^T G+\kappa I) \delta z= (1-\kappa)(G^T(\eta^\delta-\xi)),\eqno{(5.11)}
$$
where $G=(g_{k i})_{n\times 2}$ is the Jacobi matrix, and $g_{k 1}=\frac{\p u_1}{\p \alpha}(x_0, t_k)$,  $g_{k 2}=\frac{\p u_1}{\p \gamma}(x_0, t_k)$ for $k=1, 2, \cdots, n$; and
$$\eta^\delta=(u_1^{\delta}(x_0, t_1), \cdots, u_1^{\delta}(x_0, t_n))^T; \xi=(u_1[z](x_0, t_1), \cdots, u_1[z](x_0, t_n))^T.$$
By suitably choosing $\kappa\in (0, 1)$, we work out an optimal perturbation $\delta z$ by (5.11), and then we get the next iteration by linear iteration $z=z+\delta z$.\\
\indent  On the concrete inversions, we choose a Sigmoid-type function depending upon the iterations as the homotopy parameter given as
$$
\kappa(j)=\frac{1}{1+e^{\sigma\left(j-j_{0}\right)}},\eqno{(5.12)}
$$
here $j$ is the number of iterations, $j_{0}$ is the preestimated number of iterations, and $\sigma>0$ is the adjust parameter. We choose $j_{0}=5$ and $\sigma=0.9$ in all of the following computations. In addition, the forward problem is solved numerically by the finite difference scheme (5.5), and the final time is set to be $T=100$ in order to reveal the long-time behaviors of the fractional diffusion system, and the additional data are obtained at the interior point $x_{0}=0.5\in (0, 1)$. It is noted that the initial iteration is chosen as zero, i.e., $z_0=(0, 0)$ except for Ex.5.3. We refer to \cite{sunc17, zhangd13} for the detailed procedures of performance of the inversion algorithm.\\
\noindent\textbf{Example 5.1}\quad In the first numerical experiment, let $\alpha=0.8$ and $\gamma=0.25$ be the exact fractional orders, which could be suitable for some real situations where the diffusion in the immobile zone is slower than that in the mobile zone, and the exact solution of the inverse problem is expressed as $z=(0.8, 0.25)$. In addition, we take the parameters $P=5, R_{1}=R_{2}=2, \beta=0.5$, $\omega=1.5$, $\lambda=0.05$ and $\mu=0.1$ as basic settings. By substituting the exact orders into the forward problem, the solution is computed and the additional data at $x_{0}=0.5$ are obtained, with which the inversion algorithm is applied to reconstruct the fractional orders.\\
\indent The inversion results with noisy data and exact data are listed in Table 1, where $\delta$ denotes the noise level, and $\delta=0$ means that the inversion is performed with noise-free data, and $\bar{z}^{inv}:=\left(\bar{\alpha}^{inv}, \bar{\gamma}^{inv}\right)$ denotes the average inversion solution with 10-time continuous inversions, and $\bar{E}rr$ denotes the relative error in the solutions, given by $\bar{E}rr=\left\|z-\bar{z}^{inv}\right\|/\|z\|$, and $\bar{j}$ denotes the average number of iterations.
$$
\begin{array}{cccc}
	\multicolumn{4}{c} {\text { Table 1. The inversion results in Ex.5.1}} \\
	\hline \delta & \bar{z}^{i n v} & \bar{E} r r & \bar{j} \\
	\hline
	5 \%   & (0.82564596, 0.24445217) & 3.13 \mathrm{e}-2 & 20.5 \\
	1 \%   & (0.79570944, 0.25106402) & 5.27 \mathrm{e}-3 & 16.1   \\
	0.1 \% & (0.80061766, 0.25018621) & 7.69 \mathrm{e}-4 & 15.8 \\
	0.01\% & (0.79993667, 0.25000822) & 7.61 \mathrm{e}-5 & 15 \\
	0      & (0.79999999, 0.25000000) & 6.59 \mathrm{e}-10& 15   \\
	\hline
\end{array}
$$
\\
\noindent\textbf{Example 5.2} \quad In this example, we choose the model parameters as $P=1 R_{1}=R_{2}=2$, $\beta=0.5, \omega=1.5$, $\lambda=0.05$ and $\mu=0.1$, and we take $\alpha=0.75$ and $\gamma=0.75$ as the exact solution of the inverse problem, i.e., $z=(0.75, 0.75)$.  This situation could occur if the solute variations in the mobile and immobile zones have the same fractal dynamics. As done in Ex.5.1, the inversion results with noisy data and exact data are listed in Table 2.
$$
\begin{array}{cccc}
	\multicolumn{4}{c} {\text { Table 2. The inversion results in Ex.5.2}} \\
    \hline \delta & \bar{z}^{i n v} & \bar{E} r r & \bar{j} \\
	\hline
	5 \%   & (0.74747782, 0.79993481) & 4.71 \mathrm{e}-2 & 22.5 \\
	1 \%   & (0.75619961, 0.74402864) & 8.12 \mathrm{e}-3 & 20.3 \\
	0.1 \% & (0.75046744, 0.74927740) & 8.11 \mathrm{e}-4 & 18.5 \\
	0.01\% & (0.75005237, 0.74994325) & 7.28 \mathrm{e}-5 & 18   \\
	0      & (0.75000000, 0.74999999) & 1.37 \mathrm{e}-9 & 18   \\
	\hline
\end{array}
$$
\\
\textbf{Example 5.3}\quad In this example, we are concerned with a special case in which the fractional order in the immobile zone is greater than that in the mobile zone. Let $\alpha=0.3$ and $\gamma=0.8$ as the exact solution of the inverse problem, i.e., $z=(0.3, 0.8)$. The model parameters are chosen as $P=1, R_{1}=R_{2}=2$, $\beta=0.5, \omega=0.5$, $\lambda=0.05$ and $\mu=0.5$.
It is noted that the inversion results become unstable if still choosing zero as the initial iteration. The reason maybe come from the choice of the fractional orders where the order in the mobile zone is smaller than that in the immobile. However, by choosing the initial iteration as $z_0=(1, 1)$, the inversion algorithm can be realized successfully. The inversion results are listed in Table 3.
$$
\begin{array}{cccc}
	\multicolumn{4}{c} {\text { Table 3. The inversion results in Ex.5.3}} \\
	\hline \delta & \bar{z}^{i n v} & \bar{E} r r & \bar{j} \\
	\hline
	5 \%   & (0.29960169, 0.88268139) & 9.67 \mathrm{e}-2 & 28.5 \\
	1 \%   & (0.29677359, 0.81211370) & 1.46 \mathrm{e}-2 & 25.3 \\
	0.1\%  & (0.30000177, 0.80163716) & 1.91 \mathrm{e}-3 & 23.1 \\
	0.01\% & (0.30012901, 0.79979368) & 2.84 \mathrm{e}-4 & 23 \\
	0      & (0.30000000, 0.80000000) & 1.29 \mathrm{e}-10& 22\\
	\hline
\end{array}
$$
\\
\indent From Tables 1-3 it can be seen that the inversion solutions approximate to the exact solutions as the noise goes to zero, and the inversion algorithm is of numerical stability against noise in the data. The fractional orders are important to the fractal MIM solute transport model, and it could be more suitable for real situations by the inversion results that the fractional order in the mobile zone cannot be less than that in the immobile. In addition, by the natural conditions the fractional orders should be in $S_{ad}$ in theory. However, the situation could have a little change in numerical experiments. In our examples we choose $z_0=(0, 0)$ or $z_0=(1, 1)$ as the initial iteration so as to show the universality of the inversion algorithm, and the inversion results are satisfactory. Actually, if choosing $z_0=(0.1, 0.1)$ or $z_0=(0.9, 0.9)$ as the initial iteration correspondingly, the inversion results are better than those of using $z_0=(0, 0)$ or $z_0=(1, 1)$.
\section{Conclusion}
A fractal MIM solute transport model is studied from system identification. The unique existence of solution to the fractal system is discussed in mathematics by the method of Laplace transform, and the uniqueness of identifying the fractional orders is proved in the real-space of Laplace transform. Numerical inversions with noisy data are presented to demonstrate the numerical stability of the inverse problem. We will focus on the research of regularity of the solution for the forward problem, and study inverse problems of determining other parameters in the fractal system.
\section*{Appendix-Proof of Lemma 2.2}
We need to prove the convergence of the contour integral
$$
\frac{1}{2\pi i} \int_{s_{0}-i\infty}^{s_{0}+i\infty} F(s) e^{st}ds,  \eqno{(A.1)}
$$
where $F(s)$ satisfying (2.5). In the follows we denote $C$ as any positive constant if there is no specification. At first we need a convergent assertion for an infinite integral on a real-valued function, which is deduced by the comparison criterion.\\
\textbf{Lemma A.1}\quad Let $g(r)$ be a nonnegative function on $[r_0, +\infty)$ for given $r_0>0$, and be integrable on any finite interval of $[r_0, +\infty)$, and $\lim\limits_{r\rightarrow +\infty} r^p g(r)=q$. Then the integral $\int_{r_0}^{+\infty} g(r) dr$ is convergent if $p>1$ and $0\leq q<+\infty$.\\
\indent For the estimation of the contour integral (A.1), we are to utilize Cauchy integral theorem. For given angles $\theta_{1},\theta_{2}$ and a radius $\varepsilon>0$, and a infinitely large constant $R>0$, a closed curve is plotted in Figure A-1, where $L_R$ denotes the line from $s_0-iR$ to  $s_0+iR$, and $\Gamma_{R}^{+}$ denotes a finite line from $s_0+iR$ to the given point $A_1$, and $\Gamma_{+}$ denotes the line $A_1A_2$ and the circular arc $\overset\frown{A_2A_3}$, where $\theta_1\in (0, \pi/2)$ and $\theta_2\in (\pi/2, \pi)$, and there are $\Gamma_{R}^{-}$ and $\Gamma_{-}$ symmetrically corresponding to $\Gamma_R^+$ and $\Gamma_+$, respectively.
\begin{center}
\begin{tikzpicture}
\draw[-latex] (-3,0) -- (3,0) node[below] {Solid axis};   
\draw[-latex] (0,-4) -- (0,4) node[left] {Virtual axis};  
\coordinate (O) at (0,0);       
\coordinate (a) at (120:3);     
\coordinate (b) at (120:1);     
\draw[dashed] (O) -- (a) -- ++(3,0) coordinate (s0+iR) -- (-60:3) coordinate (s0-iR) -- ++(-3,0) coordinate (d);  
\draw (a) -- (b) arc (124:-124:1) -- coordinate (Gamma-) (d); 
\draw (.2,0) node[below right] {$\varepsilon$} arc (0:120:.2) node[above right] {\scriptsize$\theta_2$};             
\draw (0,.3) arc (90:120:.3) node[above] {\scriptsize$\theta_1$}; 
\node[above right] at (0.9,0) {$A_3$};   
\node[below right] at (1.5,0) {$s_0$};
\node[below right] at (s0-iR) {$s_0-\mathrm{i}R$};
\node[above right] at (s0+iR) {$s_0+\mathrm{i}R$};
\node[below left] at (b) {$A_2$};
\node[left] at (a) {$A_1$};
\node[below left] at (-0.1,0) {$O$};
\draw[->] let \p1=(a) in (.5,\y1) -- (.4,\y1) node[above] {$\Gamma_R^+$}; 
\draw[->] let \p1=(d) in (.5,\y1) -- (.6,\y1) node[below] {$\Gamma_R^-$};
\draw[->] let \p1=(s0+iR) in (\x1,.9) -- (\x1,1) node[right] {$L_R$};
\draw[->] (20:1.07) arc (20:19:1);
\draw[->] (120:2) -- (120:1.9) node[right] {$\Gamma_+$};
\draw[->] (-122:1) -- (Gamma-) node[right] {$\Gamma_-$};
\end{tikzpicture}
\begin{center}
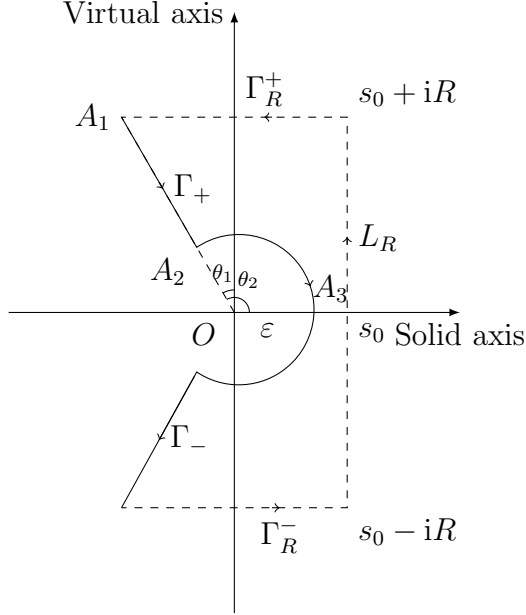
Figure A-1. A closed curve for computation of the contour integral
\end{center}
\end{center}
\vskip 0.2cm
\quad\quad From Cauchy integral theorem, it holds that
$$
\oint_{L_{R} \cup \Gamma_{R}^{+} \cup \Gamma_{+} \cup \Gamma_{-} \cup \Gamma_{R}^{-}} F(s) e^{s t} d s=0,  \eqno{(A.2)}
$$
and we get
$$
\int_{L_{R}} F(s) e^{s t} d s=-\int_{\Gamma_{R}^{\pm}} F(s) e^{s t} d s-\int_{\Gamma_{\pm}} F(s) e^{s t} d s. \eqno{(A.3)}
$$
\indent We firstly estimate the integrals on $\Gamma_{R}^{\pm}$. Denote $\Gamma_{R}^{+}=\Gamma_{R_{1}}^{+} \cup \Gamma_{R_{2}}^{+}$, where $\Gamma_{R_{1}}^{+}: s=x+i R, 0< x< s_0$; \quad $\Gamma_{R_{2}}^{+}: s=x+i R, \operatorname{Re} s \leq 0$. By the condition (2.5), there holds
$$
\begin{aligned}
	\left|\int_{\Gamma_{R_{1}}^{+}} F(s) e^{s t} d s \right|
	&\leq \int_{\Gamma_{R_{1}}^{+}} \left| F(s) \right|\cdot\left| e^{s t} \right| \left| d s \right| \\
	& \leq C \int_{0}^{s_{0}} \frac{1}{|s|} e^{x t} d x \\
	& \leq \frac{C}{R} \cdot \frac{1}{t}\left(e^{s_{0} t}-1\right).
\end{aligned}
\eqno{(A.4)}
$$
Then for any given $t>0$, there is
$$
\lim_{R \rightarrow \infty} \left|\int_{\Gamma_{R_{1}}^{+}} F(s) e^{s t} d s \right|=0.\eqno{(A.5)}
$$
\indent For the integral on $\Gamma_{R_{2}}^{+}$, we have the estimation
$$
\begin{aligned}
	\left|\int_{\Gamma_{R_{2}}^{+}} F(s) e^{s t} d s\right|
	& \leq \int_{\Gamma_{R_{2}}^{+}}\frac{C}{|s|} \cdot e^{(\mathrm{Re} s) t}|d s| \\
	& \leq \frac{C}{R} \int_{\Gamma_{R_{2}}^{+}} e^{(\mathrm{Re} s) t}|d s| \\
	& = \frac{C}{R} \int_{0}^{R \tan(\theta_{1})} e^{-x t} d x=\frac{C}{R\ t} \left(1-e^{-R \tan(\theta_{1}) t}\right).
\end{aligned}
\eqno{(A.6)}
$$
Thanks to $\theta_{1} \in (0,\frac{\pi}{2})$, there is
$$
\lim _{R \rightarrow \infty} \left|\int_{\Gamma_{R_{2}}^{+}} F(s) e^{s t} d s\right|=0.\eqno{(A.7)}
$$
Combing with (A.5) follows that
$$
\lim_{R \rightarrow \infty}\left|\int_{\Gamma_{R}^{+}} F(s) e^{s t} d s\right|=0. \eqno{(A.8)}
$$
\indent Similarly for the integral on $\Gamma_{R}^{-}$, there holds
$$
\lim_{R \rightarrow \infty}\left|\int_{\Gamma_{R}^{-}} F(s) e^{s t} d s\right|=0. \eqno{(A.9)}
$$

Next we estimate the integrals $\int_{\Gamma_{\pm}} F(s) e^{s t} d s$ in (A.3).\\
\indent As done in the above, we firstly give the estimation for the integral on $\Gamma_{+}$. Noting $\Gamma_{+}=A_1A_2 + \overset\frown{A_2A_3}$, there is
$$
|\int_{{{\Gamma}_{+}}}{F(s){{e}^{s t}}ds}|\le |\int_{A_1A_2}{F(s){{e}^{st}}}ds| + |\int_{\overset\frown{A_2A_3}}{F(s){{e}^{s t}}} ds|.
\eqno{(A.10)}
$$
For the integral $|\int_{A_1A_2} {F(s){{e}^{s t}}}ds|$, by the condition (2.5) and the polar coordinate transformation $s=r e^{i \theta_2}$ along the line $A_1A_2$, there holds
$$
\begin{aligned}
\left|\int_{A_1A_2} F(s) e^{s t} d s\right|
\leq & \int_{A_1A_2}|F(s)| e^{(\mathrm{Re} s) t}|d s| \\
	& \leq C \int_{A_1A_2} \frac{1}{|s|} e^{(\mathrm{Re} s) t}|d s| \\
	&= C \int_{\varepsilon}^{\frac{R}{\cos(\theta_{1})}} \frac{1}{r} e^{r \cos(\theta_{2}) t} d r.
\end{aligned}
\eqno{(A.11)}
$$
Noting $\theta_{2} \in (\pi/2, \pi)$, there is $\cos(\theta_{2})<0$. By utilizing Lemma A.1 where $p=2, q=0$, we deduce that for given $t>0$, the integral $\int_{\varepsilon}^{\frac{R}{\cos(\theta_{1})}} \frac{1}{r} e^{r \cos(\theta_{2}) t} d r$ is convergent as $R\rightarrow \infty$. So there exists a positive constant $C$ such that
$$
\lim_{R \rightarrow \infty} \int_{\varepsilon}^{\frac{R}{\cos(\theta_{1})}} \frac{1}{r} e^{r \cos(\theta_{2}) t} d r \leq C.\eqno{(A.12)}
$$
\indent Now we estimate the integral on the arc $\overset\frown{A_2A_3}$. There holds
$$
\left|\int_{\overset\frown{A_2A_3}} F(s) e^{s t} d s\right|
\leq C \int_{\overset\frown{A_2A_3}} \frac{1}{|s|} e^{\mathrm{Re}(s) t}|d s|.
\eqno{(A.13)}
$$
Noting that $|s|=\varepsilon$ on the circular arc, and the length of the arc is $|\overset\frown{A_2A_3}|=\frac{\theta_2 \pi}{180}\varepsilon$, we conclude that there exists a constant $C>0$ such that
$$
\int_{\overset\frown{A_2A_3}} \frac{1}{|s|} e^{\mathrm{Re}(s) t}|d s|
\leq \frac{e^{\varepsilon t}}{\varepsilon} \int_{\overset\frown{A_2A_3}} |ds|
\leq C.
\eqno{(A.14)}
$$
Therefore there exists $C>0$ such that
$$
\lim_{R \rightarrow \infty}\left|\int_{\Gamma_{+}} F(s) e^{s t} d s\right| \leq C.\eqno{(A.15)}
$$
Similarly we have
$$
\lim_{R \rightarrow \infty}\left|\int_{\Gamma_{-}} F(s) e^{s t} d s\right| \leq C.\eqno{(A.16)}
$$
\indent Based on (A.3), combing (A.15), (A.16) with  (A.8) and  (A.9), we arrive at
$$
\lim _{R \rightarrow \infty}\left|\int_{L_{R}} F(s) e^{s t} d s\right| \leq C,
\eqno{(A.17)}
$$
which means that the contour integral
$\frac{1}{2 \pi i} \int_{s_{0}-i \infty}^{s_{0}+i \infty} F(s) e^{s t} d s$ is bounded at each given $t>0$. The proof is completed.
\section*{Acknowledgements}
This work is supported by National Natural Science Foundation of China (No. 11871313), and Natural Science Foundation of Shandong Province, China (No. ZR2019MA021).


\end{document}